
\baselineskip=14pt
\parskip=10pt
\def\halmos{\hbox{\vrule height0.15cm width0.01cm\vbox{\hrule height
  0.01cm width0.2cm \vskip0.15cm \hrule height 0.01cm width0.2cm}\vrule
  height0.15cm width 0.01cm}}
 
\font\eighttt=cmtt8
\magnification=\magstephalf
\def\G{{\cal G}}

\def\C{{\cal C}}

\def\1{{\overline{1}}}
\def\2{{\overline{2}}}
\parindent=0pt
\overfullrule=0in

\def\frac#1#2{{#1 \over #2}}
\centerline
{
\bf A Very Short (Bijective!) Proof of Touchard's Catalan Identity 
}
\bigskip
\centerline
{\it By Amitai REGEV, Nathaniel SHAR, and Doron ZEILBERGER}

{\it Added March 24, 2015}: It turns out that our bijection was too nice not to have been known before.
See the insightful comments by Dominique Gouyo-Beauchamp, Kyle Petersen, and Dennis Stanton in

{\tt http://www.math.rutgers.edu/\~{}zeilberg/mamarim/mamarimhtml/touchardComm.html}
\bigskip
\hrule
\bigskip
Recall that one of the almost infinitely many definitions of the ubiquitous {\it Catalan Numbers}, $C_n$, is as the
number of elements of the set of $2n$-letter words, $w_1 \dots w_{2n}$ in the alphabet $\{-1,1\}$ that add up to zero, and 
all whose partial sums are non-negative. Let's call this set $\C_n$.

In the 1924 Toronto ICM, Jacques Touchard [T] announced (and proved) the elegant identity
$$
C_{n+1}= \sum_{ k \geq 0} {{n} \choose {2k}} 2^{n-2k} C_k \quad .
\eqno(Touchard)
$$
Here is a very short, {\it purely bijective}, proof, {\it even} nicer than Lou Shapiro's [S].

Let $f(1,1): =1$, $f(1,-1):=0$, $f(-1,1):={\bar 0}$, $f(-1,-1):= -1$,
where ${\bar 0}$ is a twin-sister of $0$, whose value in summation is also $0$. Define on $w=w_1 \dots w_{2n+2} \in \C_{n+1}$,
$$
T(w_1 \dots w_{2n+2}):=f(w_1,w_2) f(w_3,w_4) \dots f(w_{2n+1},w_{2n+2}) \quad .
$$
This is a bijection onto the set, let's call it $\G'_n$, of $n+1$-letter words in the alphabet
$\{-1,0, \bar 0, 1\}$, that sum-up to zero, have all non-negative partial sums, and {\it in addition},
the partial sum before any occurence of the letter $\bar 0$ is strictly positive.

But this latter set is in bijection with the set,  let's call it $\G_n$,  of such $n$-letter words without the last restriction,
as follows. For $w=w_1 \dots w_{n+1} \in \G'_n$, if $w_{n+1}=0$ just chop that last letter, mapping it to $w_1 \dots w_n$.
Otherwise, of course $w_{n+1}=-1$ (it can't be $1$, and it can't be ${\bar 0}$), so write $w$ as ($\bar 1 :=-1$)
$w' \, 1 \, w'' \, {\bar 1}$ (where $w' \in \G'_k$ and $w'' \in \G_{n-1-k}$ for some $0 \leq k \leq n-1$), and map it
to $w' {\bar 0} w''$.

Note that the number of elements of $\G_n$ is given by the right side of Eq. $(Touchard)$. 
Indeed, let the number of ones be $k$ ($0 \leq k \leq n/2$), then there are also $k$ minus-ones.
There are ${{n} \choose {2k}}$ ways to choose the locations of the $1$' and $-1$'s, $C_k$ ways of forming them
into a member of $\C_k$, and $2^{n-2k}$ ways of deciding which kind of zero ( $0$  or ${\bar 0}$) will occupy the
remaining $n-2k$ slots. \halmos 

{\bf Remarks}

{\bf 1.} While it is nice to give pretty bijective proofs, let us note that today, thanks to WZ proof theory, 
the {\it epistemological stature} of identities like Touchard's is the same as that of the identity $2+2=1+3$.
Indeed just copy-and-paste the line below onto a Maple session:

{\eighttt SumTools[Hypergeometric][Zeilberger](binomial(n,2*k)*2**(n-2*k)*binomial(2*k,k)/(k+1),n,k,N);}

{\bf 2.} Another way of counting $\G_n$ is to partition it according to the number of occurrences of ${\bar 0}$, say $n-k$,
then  choose the ${{n} \choose {n-k}}$ locations of the ${\bar 0}$ and `fill-in'
the remaining $k$ slots by a so-called {\it Motzkin word} of length $k$, i.e. a word in the alphabet $\{-1,0,1\}$,
whose sum is $0$, and whose partial sums are non-negative, yielding the equally elegant identity (where $M_k$ is the
number of Motzkin words of length $k$)
$$
C_{n+1}= \sum_{k=0}^{n} {{n} \choose {k}} M_k \quad .
\eqno(Motzkin)
$$
While this identity is `trivially equivalent' to quite a few known identities, 
and is `well-known to the experts', we were unable to find it in the literature.

{\bf 3.} We intentionally avoided drawing diagrams, but most human readers will probably get a better appreciation of the beauty of 
our proof by {\it drawing} a random  Dyck path in $\C_{n+1}$, and then, scanning it in  consecutive pairs,
replace $11$ (alias up-up) by an Up Step, replacing $\bar 1 \bar 1$ (alias down-down) by a Down Step,
replace $1 \bar 1$ by a green horizontal step, and replace $\bar 1 1$ by a red horizontal step.
Then $\G'_{n}$ are generalized Motzkin paths of length $n+1$, with two types of horizontal steps, green and red,
where a red horizontal step may not lie on the $x$-axis, and $\G_n$ is  the set of such $n$-step paths without this restriction.
The bijection between $\G'_{n}$ and $\G_n$ consists of removing the last step, if it is a green horizontal step
(of course it can't be a red horizontal step), and otherwise looking at the 
`Up-mate' of the last step (that is [of course] a Down step), and replacing that Up-Mate by a red horizontal step,
and at the same time deleting the above-mentioned last Down step. 

{\bf 4.} We thank Lou Shapiro for telling us that we rediscovered Touchard's identity (in its almost-equivalent
form given in Eq. $(Motzkin)$), and telling us about [S]. While we admire Shapiro's combinatorial proof, it
it is not purely bijective, and makes use of generating functions.

{\bf 5} Our bijection is a {\it renormalization-group} transformation, where we `renormalized' a word of length $2n+2$ into
a word half as long, but with more letters in the underlying alphabet. It may be interesting to see if one can
get less trivial identities by considering generalized Dyck words where the fundamental steps are drawn from a larger set 
of steps than just $\{ (1,-1), (1,1)\}$.

{\bf References}

[S] L. W. Shapiro, {\it A short proof of an identity of Touchard's concerning Catalan Numbers},
J. Combinatorial Theory (A) {\bf 20} (1976), 375-376.

[T] J. Touchard, {\it Sur certain \'equations fonctionelles}, in: Proc. Int. Math. Congress, Toronoto (1924), Vol. {\bf 1} (1928), 465-472.\hfill\break
{\tt http://www.mathunion.org/ICM/ICM1924.1/Main/icm1924.1.0465.0472.ocr.pdf} .

\vfill\eject

\hrule
\smallskip
Amitai Regev, Department of Mathematics, Weizmann Institute of Science,
amitai.regev at weizmann dot ac dot il \quad ;  \quad {\tt  http://www.wisdom.weizmann.ac.il/\~{}regev/} \quad .
\smallskip
\hrule
\smallskip
Nathaniel Shar and Doron Zeilberger, Department of Mathematics, Rutgers University (New Brunswick),
[nshar, zeilberg] at math dot rutgers dot edu \quad , \hfill\break
\quad {\tt http://www.math.rutgers.edu/\~{}zeilberg/}   \quad , \quad {\tt http://www.math.rutgers.edu/\~{}nshar/} \quad .
\smallskip
\hrule
\smallskip
{\bf ${\bf \pi}$ day, 2015}; This version: March 24, 2015.

\end